\documentclass[11pt,reqno]{amsart}
\usepackage[T1]{fontenc}
\usepackage{lmodern}
\usepackage{amsmath,amssymb,mathtools}
\usepackage{microtype}
\usepackage[a4paper,margin=29mm]{geometry}
\usepackage{xcolor}
\usepackage[numbers,sort&compress]{natbib}
\usepackage{xparse}
\usepackage[
  pagebackref=true,
  colorlinks=true,
  linkcolor=blue,
  citecolor=blue,
  urlcolor=black
]{hyperref}

\usepackage{url}
\RenewDocumentCommand{\cite}{m t* g}{%
  \IfBooleanTF{#2}{\citep[#3]{#1}}{\citep{#1}}}
\newcommand{\cites}[1]{\citep{#1}}
\renewcommand*{\backref}[1]{}
\renewcommand*{\backrefalt}[4]{%
  \ifcase #1\relax
  \or\space(cited on page~#2)%
  \else\space(cited on pages~#2)%
  \fi}

\numberwithin{equation}{section}
\newtheorem{theorem}{Theorem}[section]
\newtheorem{proposition}[theorem]{Proposition}
\newtheorem{lemma}[theorem]{Lemma}
\newtheorem{corollary}[theorem]{Corollary}
\theoremstyle{definition}

\newtheorem{example}[theorem]{Example}
\theoremstyle{remark}
\newtheorem{remark}[theorem]{Remark}
\newcommand{\ii}{\sqrt{-1}}
\newcommand{\pb}{\bar\partial}
\newcommand{\ddc}{\ii\partial\pb}
\newcommand{\F}{\mathcal F}
\newcommand{\G}{\mathcal G}
\newcommand{\HM}{\mathcal H_M}

\newcommand{\doi}[1]{\href{https://doi.org/#1}{\nolinkurl{#1}}}
\title[Pluriclosed torsion-critical metrics]{On pluriclosed torsion-critical Hermitian metrics}

\author{Dongmei Zhang}
\address{Dongmei Zhang. School of Mathematical Sciences, Chongqing Normal University, Chongqing 401331, China}
\email{{2250825921@qq.com}}

\author{Yuqin Guo}
\address{Yuqin Guo. School of Mathematical Sciences, Chongqing Normal University, Chongqing 401331, China}
\email{{1942747285@qq.com}}\thanks{Zhang is the corresponding author. Guo is supported by Chongqing graduate student research grant No. CYB25248.}

\date{}
\begin{document}
\begin{abstract}
  In our previous work, we studied the squared $L^2$-norm functionals $\mathcal G$ of Gauduchon's torsion $1$-form and $\mathcal F$ of the Chern torsion over the space of all Hermitian metrics with unit volume on a compact complex manifold. We introduced torsion-critical metrics as critical points of $\mathcal F$.
  
  In this paper, We prove that a pluriclosed and astheno-K\"ahler metric is torsion-critical if and only if it is K\"ahler.
  In particular, every compact pluriclosed torsion-critical Hermitian threefold is K\"ahler.
  We also give explicit examples and derive geometric consequences under additional curvature assumptions.
  \end{abstract}

\maketitle

\tableofcontents

\section{Introduction and statement of results}

Let $(M^n,g)$ be a compact Hermitian manifold of complex dimension $n$.
Denote by $\omega$ its K\"ahler form, $\nabla$ its Chern connection, and $T$ the torsion of $\nabla$. 
The Gauduchon torsion $1$-form $\eta $, defined as the trace of the torsion tensor T, is the global \((1,0)\)-form \(\eta\) on \(M^n\) given by
\begin{equation}\label{eq:eta-intro}
  \partial\omega^{n-1}=-\eta\wedge\omega^{n-1}.
\end{equation}
For $n\geq2$, the condition $\eta=0$ is equivalent to balancedness, that is, $d\omega^{n-1}=0$. If the torsion T vanishes, then the metric g is Kähler, that is, $d\omega=0$.
Gauduchon initiated the variational study of the torsion $1$-form in \cite{Gauduchon84}. 
He considered the functional
\begin{equation}\label{eq:functionals}
  \G(g)=V^{\frac{1-n}{n}}\int_M|\eta|^2\,dv,
  \qquad g\in\HM,
\end{equation}
he derived its Euler--Lagrange equation, and proved that its critical metrics on compact complex surfaces are balanced. 

In a previous work\cite{ZZ23}*{Proposition~1} ,we extended this conclusion to all complex dimensions $n\geq2$. 
In the same work, we introduced \emph{torsion-critical metrics} which are Hermitian metrics on a compact complex manifold $M^n$ taht are critical points of the functional of $L^2$-norm of the Chern torsion tensor $T$:
\begin{equation}\label{eq:F-intro}
  \F(g)=V^{\frac{1-n}{n}}\int_M|T|^2\,dv,
  \qquad g\in\HM.
\end{equation}
Throughout this paper, we use their normalization$|T|^2=\sum_{i,j,k}|T^j_{ik}|^2.$
Where ${\mathcal H}_M$ is the set of all Hermitian metrics on the complex manifold $M^n$, $dv=\frac{\omega^n}{n!}$ the volume form, and $V=\int_Mdv$ the volume. Note that the factor in front of the integral is simply to make the functional independent of the scaling of metric (by constant multiples). Equivalently, one can drop this $V$ factor but restrict to the subset of all Hermitian metrics with unit volume.

One may also consider the squared norm of the Bismut torsion. Let $T^{\mathrm B}$ denote the torsion tensor of the Bismut connection. 
The Bismut connection is the unique Hermitian connection whose torsion is totally skew-symmetric \cites{Bismut,Gauduchon97},  for flatness questions involving Hermitian connections, see \cites{VYZ,WYZ}. 
Its torsion three-form is $d^c \omega $ up to sign, where $d^c=\sqrt{-1}(\overline{\partial}-\partial) $; the sign depends on the convention for the fundamental form. 
Hence, with the normalization of the Chern torsion used in this paper,
$$
 |T^{\mathrm B}|^2=6|d\omega|^2=12|\partial\omega|^2=6|T|^2.
$$
It follows that
$$
 V^{\frac{1-n}{n}} \int_M |T^{\mathrm B}|^2\,dv =6\mathcal F(g)
$$
for every $g\in\mathcal H_M$.
Thus the two functionals have the same critical points with respect to all Hermitian variations.

The critical metrics of $\mathcal G$ are all balanced, and hence all its critical values are zero. In contrast, $\mathcal F$ admits non-Kähler critical metrics in addition to Kähler critical metrics, at which its value is zero.
Zhang and Zheng constructed non-K\"ahler examples and proved several non-existence results in \cite{ZZ23}. 
Their subsequent work \cite{ZZ25} addresses the compact Chern flat case: the complex Lie group arising as the universal cover of a non-K\"ahler Chern flat torsion-critical manifold must be semisimple. Conversely, canonical metrics on semisimple complex Lie groups yield torsion-critical metrics on compact quotients whenever such quotients exist. 
Thus, torsion-criticality alone does not imply the K\"ahler condition. 
Related variational problems within conformal classes were studied by Angella, Istrati, Otiman, and Tardini \cite{AIOT}, while Kawamura \cite{Kawamura} considered a torsion variational problem in the almost Hermitian setting.
In this paper, we fix an integrable complex structure and consider general variations of compatible Hermitian metrics.

One of the questions raised in \cite{ZZ23} concerns
torsion-critical metrics satisfying the \emph{pluriclosed}
condition
\begin{equation}\label{eq:SKT-intro}
  \partial\pb\omega=0.
\end{equation}
Such metrics are also called strong K\"ahler with torsion, or SKT.
They arise naturally in the theory of Hermitian connections with torsion; see \cites{Bismut,Gauduchon97,Strominger} and the survey \cite{FTsurvey}. 
Non-K\"ahler pluriclosed metrics on complex three-dimensional nilmanifolds were studied in \cite{FPS}. Our three-dimensional result implies that such metrics cannot be torsion-critical.
The relation with symplectic forms taming the complex structure was investigated in \cite{EFV}.
Other restrictions on the existence of non-K\"ahler pluriclosed metrics appear in \cite{Chiose}. The pluriclosed condition is also preserved by the geometric flow introduced by Streets and Tian \cite{StreetsTian}.

We study this question under the additional \emph{astheno-K\"ahler} condition
\begin{equation}\label{eq:AK-intro}
  \partial\pb\omega^{n-2}=0,
  \qquad n\geq3.
\end{equation}
Introduced by Jost and Yau \cite{JostYau}, this condition has been studied from the perspectives of both construction and rigidity; see, for example, \cites{FTastheno,FGV}. Related conditions on powers of the fundamental form occur in the study of generalized Gauduchon metrics \cite{FinoUgarte}.
The pluriclosed and astheno-K\"ahler conditions coincide in complex dimension three but are distinct in higher dimensions. In this paper, simultaneous conditions are imposed on the same metric.

Our argument uses balancedness as an intermediate step. Balanced metrics were introduced systematically by Michelsohn \cite{Michelsohn}, and their interaction with other torsion conditions is known to be an important source of Kähler rigidity; see \cites{AlexIvanov,FV}. 
Note, however, that the existence of a single metric satisfying two conditions is a strictly stronger requirement than the existence of a manifold admitting different metrics that satisfy them separately. 
In fact, \cites{FGV,LatorreUgarte} give examples of non-Kähler compact complex manifolds that admit both a balanced and an astheno-Kähler metric.

Our main result gives a partial answer to the pluriclosed problem posed in \cite{ZZ23}.

\begin{theorem}\label{thm:main}
Let $(M^n,g)$ be a compact Hermitian manifold, with $n\geq3$. Suppose that $g$ is both pluriclosed and astheno-K\"ahler. Then $g$ is torsion-critical if and only if it is K\"ahler.
\end{theorem}

In particular, in the complex three-dimensional case, we obtain the following consequence:

\begin{corollary}\label{cor:three}
On a compact complex threefold, a Hermitian metric is pluriclosed and torsion-critical if and only if it is K\"ahler.
\end{corollary}

In complex dimension two, the identity $\F=2\G$ under the above normalization and Gauduchon's theorem imply that every torsion-critical metric is K\"ahler. 
In complex dimension one, every Hermitian metric is K\"ahler. Thus, the characterization in Corollary~\ref{cor:three} holds in all complex dimensions at most three. 
In contrast to the Chern-flat results of \cite{ZZ25}, Theorem~\ref{thm:main} requires no curvature or homogeneity assumptions.

The proof rests on the identity
\begin{equation}\label{eq:intro-identity}
  \F(g)-2\G(g)
  =
  \frac{2V^{\frac{1-n}{n}}}{(n-2)!}
  \int_M\ddc\omega\wedge\omega^{n-2},
\end{equation}
which holds for every Hermitian metric when $n\geq3$. For a smooth curve $g_t\in\HM$, write
$$
  D\F_g(\gamma)=\left.\frac{d}{dt}\right|_{t=0}\F(g_t),
$$
and define $D\G_g(\gamma)$ similarly. Here $\gamma$ is an arbitrary smooth real $(1,1)$-form. At a pluriclosed metric, differentiating \eqref{eq:intro-identity} and integrating by parts gives
\begin{equation}\label{eq:intro-variation}
  (D\F_g-2D\G_g)(\gamma)
  =
  \frac{2V^{\frac{1-n}{n}}}{(n-2)!}
  \int_M\gamma\wedge\ddc\omega^{n-2}.
\end{equation}
The astheno-K\"ahler condition therefore yields $D\F_g=2D\G_g$. If $g$ is torsion-critical, it is consequently
a critical point of $\G$ and hence balanced by \cite{ZZ23}*{Proposition~1}. Equation~\eqref{eq:intro-identity} then implies $\F(g)=0$, so $g$ is K\"ahler. The converse follows because every K\"ahler metric minimizes the nonnegative functional $\F$.

For a pluriclosed metric, the identity $\mathcal F=2\mathcal G$ does not by itself imply $D\mathcal F=2D\mathcal G$.
The astheno-K\"ahler condition ensures equality of these first variations. Without this additional condition, the pluriclosed torsion-critical problem remains unresolved by our argument in complex dimensions $n\geq4$.

The paper is organized as follows. In Section $2$, we prove Theorem~\ref{thm:main} and derive several consequences.
In Section $3$, we give some examples illustrating the theorem, and discuss the role of the astheno-K\"ahler condition.
In Section $4$, we apply the main theorem to obtain geometric consequences under Chern flatness or positive Chern holomorphic bisectional curvature.


\section{Proof of Theorem 1.1}\label{sec:identities}

Throughout this section, let $(M^n,g)$ be a compact Hermitian manifold. We use the same notation $g$ for its complex bilinear extension to the complexified tangent bundle. Let $\nabla$ be the Chern connection, with torsion defined by
$$
  T(X,Y)=\nabla_XY-\nabla_YX-[X,Y].
$$
The Chern torsion has no mixed component:
$T(X,\overline{Y})=0$ for $X,Y\in T^{1,0}M$.

In local holomorphic coordinates $(z_1,\ldots,z_n)$, set
$$
  \partial_i=\frac{\partial}{\partial z_i},
  \qquad
  g_{i\bar j}=g(\partial_i,\overline{\partial_j}),
$$
and write $(g^{\bar j i})$ for the inverse of $(g_{i\bar j})$,
so that
$$
  \sum_j g_{i\bar j}g^{\bar j k}=\delta_i^k.
$$
We adopt the convention
$$
  \nabla_{\partial_k}\partial_i
  =\sum_j\Gamma^j_{ik}\partial_j,
  \qquad
  \Gamma^j_{ik}
  =\sum_{\ell}g^{\bar\ell j}\partial_k g_{i\bar\ell}.
$$
Consequently, the coefficients in
$T(\partial_i,\partial_k)=\sum_jT^j_{ik}\partial_j$ satisfy
\begin{equation}\label{eq:T}
  T^j_{ik}
  =\Gamma^j_{ki}-\Gamma^j_{ik}
  =\sum_{\ell}g^{\bar\ell j}
    \bigl(\partial_i g_{k\bar\ell}
          -\partial_k g_{i\bar\ell}\bigr).
\end{equation}

The fundamental form of $g$ is
$$
  \omega=\sqrt{-1}\sum_{i,j}
  g_{i\bar j}\,dz_i\wedge d\bar z_j.
$$
With the sign convention in \eqref{eq:eta-intro},
Gauduchon's torsion $1$-form is
$\eta=\sum_i\eta_i\,dz_i$, where
\begin{equation}\label{eq:eta}
  \eta_i
  =\sum_kT^k_{ki}
  =\sum_k\bigl(\Gamma^k_{ik}-\Gamma^k_{ki}\bigr).
\end{equation}

For pointwise computations, choose a local unitary frame
$\{e_i\}_{i=1}^n$ of $T^{1,0}M$ with dual coframe
$\{\varphi_i\}_{i=1}^n$. We use the same component notation
relative to the chosen frame:
\begin{equation}\label{eq:conventions}
  \omega=\sqrt{-1}\sum_i
  \varphi_i\wedge\overline{\varphi_i},
  \qquad
  T(e_i,e_k)=\sum_jT^j_{ik}e_j,
  \qquad
  \eta=\sum_i\left(\sum_kT^k_{ki}\right)\varphi_i.
\end{equation}

We use the complex orientation and set $dv=\omega^n/n!$.
Norms of differential forms are induced by $g$, and the
Hodge star is extended complex linearly. For the Chern
torsion, we use the normalization
$$
 |T|^2=\sum_{i,k,j}|T^j_{ik}|^2.
$$

Since $T^j_{ik}=-T^j_{ki}$, the Chern structure equations yield
\begin{equation}\label{eq:Tnorm}
  \partial\omega
  =\sqrt{-1}\sum_{i<k}\sum_j
    T^j_{ik}\,
    \varphi_i\wedge\varphi_k\wedge\overline{\varphi_j},
  \qquad
  |\partial\omega|^2
  =\sum_{i<k}\sum_j|T^j_{ik}|^2
  =\frac12|T|^2.
\end{equation}

\begin{lemma}\label{lem:pointwise}
On a Hermitian manifold of complex dimension $n\geq3$,
\begin{equation}\label{eq:pointwise}
 \frac{\ii\partial\omega\wedge\pb\omega\wedge\omega^{n-3}}{(n-3)!}
 =\left(|\eta|^2-|\partial\omega|^2\right)dv.
\end{equation}
\end{lemma}

\begin{proof}
The argument is pointwise. Put $\alpha=\partial\omega$. On a Hermitian vector space, the Lefschetz decomposition gives
\begin{equation}\label{eq:lef}
 \alpha=\alpha_0+\beta\wedge\omega,
 \qquad \alpha_0\wedge\omega^{n-2}=0,
\end{equation}
where $\alpha_0$ is a primitive $(2,1)$-form and $\beta$ is a $(1,0)$-form. First take $\partial$ of $\omega^{n-1}$. Since $\omega$ has total degree 2, no minus sign appears in the product rule:
$$
\partial \omega^{n-1} = (n-1)\partial \omega \wedge \omega^{n-2}.
$$By \eqref{eq:eta-intro},
$$
 (n-1)\alpha\wedge\omega^{n-2}=-\eta\wedge\omega^{n-1}.
$$
Wedging $1$-forms with $\omega^{n-1}$ is injective. Substitution of \eqref{eq:lef} therefore gives $\beta=-\eta/(n-1)$. Since the Lefschetz decomposition is orthogonal and $|\beta\wedge\omega|^2=(n-1)|\beta|^2$, we obtain
\begin{equation}\label{eq:orthogonal}
 |\alpha|^2=|\alpha_0|^2+\frac{|\eta|^2}{n-1}.
\end{equation}

The Hodge star operator $*$ is the complex linear extension of the real Hodge star operator, and therefore
$$
u \wedge *\overline{u} = |u|^2 dv.
$$
For clarity, with our conventions the Hodge-star formula for a primitive $(p,r)$-form $u$ of degree $k=p+r$ is
$$
 *u=(-1)^{k(k+1)/2}(\ii)^{p-r}
 \frac{u\wedge\omega^{n-k}}{(n-k)!}.
$$
Applying this to $\overline{\alpha_0}$, and using $\alpha_0\wedge *\overline{\alpha_0}=|\alpha_0|^2dv$, yields
\begin{equation}\label{eq:primitive}
 \frac{\ii\alpha_0\wedge\overline{\alpha_0}\wedge\omega^{n-3}}{(n-3)!}
 =-|\alpha_0|^2dv.
\end{equation}
Also, a direct calculation in a unitary coframe gives
$$
 \frac{\ii\eta\wedge\bar\eta\wedge\omega^{n-1}}{(n-1)!}
 =|\eta|^2dv.
$$

Since $\omega$ is real, the primitivity of $\alpha_0$ gives
$$
 \alpha_0\wedge\omega^{n-2}=0,
 \qquad
 \overline{\alpha_0}\wedge\omega^{n-2}=0.
$$
Thus the mixed terms in the expansion of
$\alpha\wedge\overline{\alpha}\wedge\omega^{n-3}$ vanish,
and hence
$$
 \alpha\wedge\overline{\alpha}\wedge\omega^{n-3}
 =
 \alpha_0\wedge\overline{\alpha_0}\wedge\omega^{n-3}
 +\beta\wedge\overline{\beta}\wedge\omega^{n-1}.
$$

For the third term, take the complex conjugate of \eqref{eq:lef} and use $\bar{\omega}=\omega$, giving
$
\overline{\alpha_0} \wedge \omega^{n-2} = 0.
$
Therefore
$$
\beta \wedge \overline{\alpha_0} \wedge \omega^{n-2} = 0.
$$

Consequently,
\begin{align*}
 \frac{\ii\alpha\wedge\bar\alpha\wedge\omega^{n-3}}{(n-3)!}
 &=\frac{\sqrt{-1} \alpha_0 \wedge \overline{\alpha_0} \wedge \omega^{n-3}}{(n-3)!}+\frac{\sqrt{-1} \beta \wedge \bar{\beta} \wedge \omega^{n-1}}{(n-3)!}  \\
 &=\left(-|\alpha_0|^2+\frac{(n-1)(n-2)}{(n-1)^2}|\eta|^2\right)dv\\
 &=\left(|\eta|^2-|\alpha|^2\right)dv,
\end{align*}
where the last equality follows from \eqref{eq:orthogonal}. Since $\bar\alpha=\pb\omega$, this is \eqref{eq:pointwise}.
\end{proof}

\begin{proposition}\label{prop:integral}
Let $(M^n,g)$ be a compact Hermitian manifold, with $n\geq3$. Then
\begin{equation}\label{eq:integral}
 \int_M\ddc\omega\wedge\omega^{n-2}
 =(n-2)!\int_M\left(\frac12|T|^2-|\eta|^2\right)dv.
\end{equation}
In particular, \eqref{eq:intro-identity} holds, and every pluriclosed metric satisfies $\F(g)=2\G(g)$.
\end{proposition}

\begin{proof}
The form $\pb\omega$ has odd total degree. The product rule, followed by interchange of the two three-forms, gives
$$
 \partial(\pb\omega\wedge\omega^{n-2})
 =\partial\pb\omega\wedge\omega^{n-2}
 +(n-2)\partial\omega\wedge\pb\omega\wedge\omega^{n-3}.
$$
The form $\bar\partial\omega\wedge\omega^{n-2}$ has type $(n-1,n)$, so its $\bar\partial$-derivative vanishes. Its $\partial$-derivative therefore equals its exterior derivative, whose integral is zero by Stokes' theorem.
Multiplying by $\ii$ and applying Lemma~\ref{lem:pointwise}, we obtain
\begin{align*}
 \int_M\ddc\omega\wedge\omega^{n-2} 
 &=-(n-2)\int_M\sqrt{-1}\partial\omega\wedge\pb\omega\wedge\omega^{n-3} \\
 &=(n-2)!\int_M\bigl(|\partial\omega|^2-|\eta|^2\bigr)dv.
\end{align*}
Equation \eqref{eq:Tnorm} proves \eqref{eq:integral}. Multiplication by $2V^{\frac{1-n}{n} }/(n-2)!$ proves \eqref{eq:intro-identity}. Its right-hand side vanishes if $g$ is pluriclosed.
\end{proof}

\begin{remark}\label{rem:two}
  In complex dimension two, \eqref{eq:eta-intro} reduces to
$$
    \partial\omega=-\eta\wedge\omega.
$$
  Since the map $\beta\mapsto\beta\wedge\omega$ is an isometry
  on $(1,0)$-forms in this dimension, \eqref{eq:Tnorm} gives
 $$
    |T|^2=2|\partial\omega|^2=2|\eta|^2.
$$
  Thus $\F=2\G$ on the entire space $\HM$, without any
  pluriclosedness assumption.
\end{remark}

Let $\gamma$ be an arbitrary smooth real $(1,1)$-form and set $\omega_t=\omega+t\gamma$. By compactness, $\omega_t$ is positive for sufficiently small $|t|$ and hence determines a smooth family of Hermitian metrics $g_t\in\HM$. 
No differential constraint is imposed on $\gamma$. Since the complex structure is fixed,
$$
  \left.\frac{d}{dt}\right|_{t=0}
  \partial\bar\partial\omega_t
  =\partial\bar\partial\gamma,
  \qquad
  \left.\frac{d}{dt}\right|_{t=0}\omega_t^{n-2}
  =(n-2)\gamma\wedge\omega^{n-3}.
$$

\begin{lemma}\label{lem:IBP}
For a real $(1,1)$-form $\gamma$ and a smooth $(n-2,n-2)$-form $\Omega$ on $M$,
\begin{equation}\label{eq:IBP}
 \int_M\ii\partial\pb\gamma\wedge\Omega
 =\int_M\gamma\wedge\ii\partial\pb\Omega.
\end{equation}
\end{lemma}

\begin{proof}
Two applications of Stokes' theorem give
\begin{align*}
 0&=\int_M\partial(\pb\gamma\wedge\Omega)
 =\int_M\partial\pb\gamma\wedge\Omega
   -\int_M\pb\gamma\wedge\partial\Omega,\\
 0&=\int_M\pb(\gamma\wedge\partial\Omega)
 =\int_M\pb\gamma\wedge\partial\Omega
   -\int_M\gamma\wedge\partial\pb\Omega.
\end{align*}
In the second line we used $\pb\partial=-\partial\pb$. Combining the two equations and multiplying by $\ii$ proves the formula.
\end{proof}

\begin{proposition}\label{prop:variation}
For every $g\in\HM$ of complex dimension $n\geq3$, one has
\begin{align}\label{eq:full-variation1}
 (D\F_g-2D\G_g)(\gamma)
 =\frac{2V^{\frac{1-n}{n} }}{(n-2)!}\bigg\{&\int_M\gamma\wedge\ddc\omega^{n-2}\notag\\
 &+(n-2)\int_M\ddc\omega\wedge\gamma\wedge\omega^{n-3}\notag\\
 &+\frac{\frac{1-n}{n} \int_M\ddc\omega\wedge\omega^{n-2}}{V}\int_M\frac{\gamma\wedge\omega^{n-1}}{(n-1)!}\bigg\}.
\end{align}
If $g$ is pluriclosed, this reduces to
\begin{equation}\label{eq:variation}
 (D\F_g-2D\G_g)(\gamma)
 =\frac{2V^{\frac{1-n}{n} }}{(n-2)!}\int_M\gamma\wedge\ddc\omega^{n-2}.
\end{equation}
\end{proposition}

\begin{proof}
  Let $\omega_t=\omega+t\gamma$, and let
$$
   V_t=\int_M\frac{\omega_t^n}{n!}.
$$
  Since the complex structure is fixed, differentiation gives
  \begin{align*}
   \left.\frac{d}{dt}\right|_{t=0}
   \int_M \ii\partial\pb\omega_t\wedge\omega_t^{n-2}
   &=
   \int_M \ii\partial\pb\gamma\wedge\omega^{n-2}\\
   &\quad +(n-2)\int_M
   \ii\partial\pb\omega\wedge\gamma\wedge\omega^{n-3}.
  \end{align*}
  By Lemma~\ref{lem:IBP}, the first integral on the right is
$\int_M\gamma\wedge\ii\partial\pb\omega^{n-2}.$
  Moreover,
$$
   \left.\frac{d}{dt}\right|_{t=0}V_t
   =
   \int_M\frac{\gamma\wedge\omega^{n-1}}{(n-1)!}.
$$
  Differentiating \eqref{eq:intro-identity} and applying
  the product rule now gives \eqref{eq:full-variation1}.
  
  If $g$ is pluriclosed, then $\partial\pb\omega=0$.
  The second and third terms in
  \eqref{eq:full-variation1} vanish, and
  \eqref{eq:variation} follows.
  \end{proof}

\begin{proposition}\label{prop:equivalence}
Let $g$ be a pluriclosed metric on a compact complex manifold of dimension $n\geq3$. The following statements are equivalent:
\begin{enumerate}
\item $g$ is astheno-K\"ahler;
\item $D\F_g=2D\G_g$ as linear functionals on all smooth real $(1,1)$-forms;
\item $g$ is a critical point of $\F-2\G$ on $\HM$.
\end{enumerate}
\end{proposition}

\begin{proof}
The implication (1)$\Rightarrow$(2) follows from \eqref{eq:variation}, and (2) and (3) are equivalent by definition. 
To prove $(2)\Rightarrow(1)$, set
$$\Psi=\ii\partial\pb\omega^{n-2}.$$
Since $\omega$ is real and $\pb\partial=-\partial\pb$, the form $\Psi$ is real and of type $(n-1,n-1)$.
Its Hodge dual $*\Psi$ is a smooth real $(1,1)$-form and is an admissible variation direction. Under condition $(2)$, equation \eqref{eq:variation}, applied to $\gamma=*\Psi$, gives
$$
 0=\int_M *\Psi\wedge\Psi=\int_M|\Psi|^2dv.
$$
Here the middle equality follows because both forms have even degree. Hence $\Psi=0$, so $g$ is astheno-K\"ahler.
\end{proof}

\begin{remark}
The equivalence concerns equality of the full differentials. A variation that remains pluriclosed satisfies $\partial\pb\gamma=0$ to first order. Equality of derivatives only for such directions is weaker and cannot be used to conclude $D\G_g=0$ on $\HM$. Proposition~\ref{prop:variation} avoids this issue by computing at an arbitrary base metric first.
\end{remark}

We use the following theorem, and follow the conventions established above.

\begin{theorem}[Zhang--Zheng {\cite{ZZ23}*{Proposition 1}}]\label{thm:G}
On a compact complex manifold of complex dimension $n\geq2$, a Hermitian metric is a critical point of $\G$ on $\HM$ if and only if it is balanced.
\end{theorem}

The hard part of the theorem is proving that critical metrics are balanced. Conversely, a balanced metric satisfies $\mathcal G(g)=0$ and hence is an absolute minimum of the nonnegative functional $\mathcal G$.
To apply this theorem, we need $D\mathcal G_g(\gamma)=0$ for every smooth real $(1,1)$-form $\gamma$. This is why \eqref{eq:variation} is established for arbitrary Hermitian variations.

\begin{proof}[Proof of Theorem~\ref{thm:main}]
Suppose first that $g$ is torsion-critical, pluriclosed, and astheno-K\"ahler. For any smooth real $(1,1)$-form $\gamma$, Proposition~\ref{prop:variation} gives
$$
 0=D\F_g(\gamma)=2D\G_g(\gamma).
$$
Hence $g$ is a critical point of $\G$ on $\HM$. Theorem~\ref{thm:G} implies $\eta=0$. In particular, $\G(g)=0$. By Proposition~\ref{prop:integral} and pluriclosedness,
$$
 \F(g)=2\G(g)=0.
$$
The function $|T|^2$ is smooth and nonnegative, and $V>0$, so $T$ vanishes identically. Equation \eqref{eq:Tnorm} gives $\partial\omega=0$; since $\omega$ is real, also $\pb\omega=0$. Thus $d\omega=0$ and $g$ is K\"ahler.

Conversely, a K\"ahler metric satisfies both differential conditions. It also satisfies $\F(g)=0$ and therefore is an absolute minimum of the nonnegative smooth functional $\F$. Its differential vanishes in every direction, so it is torsion-critical.
\end{proof}

\begin{proof}[Proof of Corollary~\ref{cor:three}]
When $n=3$, one has $\omega^{n-2}=\omega$. Thus pluriclosedness and the astheno-K\"ahler condition coincide, and Theorem~\ref{thm:main} applies.
\end{proof}

\begin{corollary}\label{cor:low}
  Every compact pluriclosed torsion-critical Hermitian manifold of complex dimension at most three is Kähler.
  In complex dimension two, the pluriclosed assumption is unnecessary.
  \end{corollary}
  
  \begin{proof}
  The case $n=3$ follows from Corollary \ref{cor:three}. For $n=2$, Remark \ref{rem:two} gives $\mathcal F=2\mathcal G$,
  so Theorem \ref{thm:G} implies that every torsion-critical metric is balanced, hence Kähler.
  The case $n=1$ is automatic.
  \end{proof}

We also recover the following well-known consequence.

\begin{corollary}\label{cor:balanced}
  Let $(M^n,g)$ be a compact Hermitian manifold. If $g$ is both balanced and pluriclosed, then it is K\"ahler.
\end{corollary}

\begin{proof}
For $n\geq3$, balancedness gives $\G(g)=0$, and Proposition~\ref{prop:integral} gives $\F(g)=2\G(g)=0$. Thus $T=0$. Dimensions one and two are immediate as in Corollary~\ref{cor:low}.
\end{proof}

Let $\nabla^{\mathrm B}$ denote the Bismut connection (or Strominger connection), we call $g$ \emph{Strominger torsion parallel} (STP) if
$\nabla^{\mathrm B}T=0,$
where $T$ is the Chern torsion; see \cite{ZZ23}.
This convention also appears in work relating the Bismut connection to pluriclosed geometry; see \cite{ZhaoZheng}.

\begin{corollary}\label{cor:STP}
  Let $(M^n,g)$ be a compact Hermitian manifold. If $g$ is pluriclosed, torsion-critical, and STP, then it is K\"ahler.
\end{corollary}

\begin{proof}
By \cite{ZZ23}*{Proposition 3}, every compact torsion-critical STP metric is balanced, Corollary \ref{cor:balanced} implies that it is Kähler.
\end{proof}

\begin{remark}
The argument for Corollary~\ref{cor:STP} uses the STP theorem of \cite{ZZ23} rather than Theorem~\ref{thm:main}; it does not assume that every pluriclosed STP metric is astheno-K\"ahler. Likewise, the non-K\"ahler compact Chern-flat critical examples in \cite{ZZ25} do not contradict the main theorem. Compact Chern-flat metrics are balanced, as recalled in that work using \cite{YangZheng}, so Corollary~\ref{cor:balanced} already excludes pluriclosedness for those non-K\"ahler metrics.
\end{remark}


\section{Examples and the astheno-K\"ahler condition}\label{sec:scope}

As shown in the introduction, the $L^2$-norm functional of the Bismut torsion equals $6\mathcal F$ and therefore has the same critical metrics. We now give some examples and discuss the role of the astheno-Kähler condition.

\begin{example}[A flat complex three-torus]\label{ex:torus}
Let $M=\mathbb C^3/(\mathbb Z+\ii\mathbb Z)^3$ and write $z_j=x_j+\ii y_j$. The Euclidean metric and its fundamental form are
$$
 g_0=\sum_{j=1}^3(dx_j^2+dy_j^2),\qquad
 \omega_0=\frac{\ii}{2}\sum_{j=1}^3dz_j\wedge d\bar z_j.
$$
Both descend to $M$, and $d\omega_0=0$. Thus $g_0$ is Kähler and hence pluriclosed. Since $\mathcal F(g_0)=0$ and $\mathcal F\geq 0$, $g_0$ is torsion-critical. This is an example of Corollary~\ref{cor:three}.
\end{example}

\begin{example}[Complex projective space]\label{ex:projective}
For every $n\geq4$, consider $\mathbb{CP}^n$ with the Fubini--Study form, written in an affine chart as
$$
 \omega_{\mathrm{FS}}=\ii\partial\pb
 \log\left(1+\sum_{j=1}^n|z_j|^2\right).
$$
These local expressions define a global positive closed real $(1,1)$-form. Thus $\partial\pb\omega_{\mathrm{FS}}=0$ and $\partial\pb\omega_{\mathrm{FS}}^{n-2}=0$. The Chern torsion vanishes, so the metric is an absolute minimum of $\F$ on $\HM$. It therefore satisfies all the hypotheses of Theorem~\ref{thm:main}. This also shows that the theorem does not require flatness: the Fubini--Study metric has positive holomorphic bisectional curvature.
\end{example}

The next example shows why torsion-criticality cannot be omitted. We use the usual description of a primary Kodaira surface as a nilmanifold; see \cites{Hasegawa,Salamon} for the nilmanifold setting.

\begin{example}[A primary Kodaira surface times a flat complex torus]\label{ex:kodaira}
  For every $n\geq 3$, we construct a non-Kähler Hermitian metric on a compact complex $n$-manifold that is both pluriclosed and astheno-Kähler. We also verify directly that it is not torsion-critical.

To construct it, let $G=\mathbb R^4$ with multiplication
$$
 (x,y,z,t)(x',y',z',t')
 =(x+x',y+y',z+z'+xy',t+t').
$$
The integer points form a discrete cocompact subgroup $\Gamma$. Indeed, left multiplication by integer points reduces $x,y,t$ to bounded fundamental intervals and then does the same for $z$. On $S=\Gamma\backslash G$, the left-invariant forms
$$
 e^1=dx,\quad e^2=dy,\quad e^3=dz-x\,dy,\quad e^4=dt
$$
are globally defined and satisfy
$$
 de^1=de^2=de^4=0,\qquad de^3=-e^1\wedge e^2.
$$
Let $\varphi^1=e^1+\ii e^2$ and $\varphi^2=e^3+\ii e^4$ be of type $(1,0)$. Then
\begin{equation}\label{eq:kodaira-structure}
 d\varphi^1=0,\qquad
 d\varphi^2=-\frac{\ii}{2}\varphi^1\wedge\overline{\varphi^1}.
\end{equation}
There are no $(0,2)$ terms, so the complex structure is integrable. Take
$$
 g_S=\sum_{a=1}^4 e^a\otimes e^a,\qquad
 \alpha=\frac{\ii}{2}\left(\varphi^1\wedge\overline{\varphi^1}
             +\varphi^2\wedge\overline{\varphi^2}\right).
$$
Since $\alpha=e^1\wedge e^2+e^3\wedge e^4$, we have
$$
 d\alpha=-e^1\wedge e^2\wedge e^4\neq0.
$$

By (\ref{eq:kodaira-structure}), $d\varphi^2$ is a real $(1,1)$-form. Hence
$$
\bar\partial\varphi^2=\partial\overline{\varphi^2}=-\frac{\sqrt{-1}}{2} \varphi^1\wedge\overline{\varphi^1},
\qquad
\partial\varphi^2 =\bar\partial\overline{\varphi^2}=0.
$$
Since $d\varphi^1=0$, it follows that
$$
\bar\partial\alpha =\frac14\varphi^1\wedge\overline{\varphi^1} \wedge\overline{\varphi^2},
\qquad
\partial\bar\partial\alpha =-\frac{\sqrt{-1}}8 \varphi^1\wedge\overline{\varphi^1} \wedge\varphi^1\wedge\overline{\varphi^1} =0.
$$
Thus $\alpha$ is pluriclosed.

Let $X=\mathbb C^{n-2}/(\mathbb Z+\sqrt{-1}\mathbb Z)^{n-2}$ with its flat Kähler form $\beta$. On $M=S\times X$, set
$$
\omega=\pi_S^*\alpha+\pi_X^*\beta,
$$
where $\pi_S$ and $\pi_X$ are the projections onto the two factors. Since these projections are holomorphic,
$$
\partial\bar\partial\omega
=\pi_S^*(\partial\bar\partial\alpha)
+\pi_X^*(\partial\bar\partial\beta)=0.
$$
Moreover,
$$
d\omega=\pi_S^*(d\alpha)\neq 0,
$$
so the product metric is not Kähler.

Since $S$ has complex dimension two, $\alpha^3=0$ and $d(\alpha^2)=0$. For $k\geq 2$, we therefore have
$$
\omega^k
=\pi_X^*(\beta^k)
+k\,\pi_S^*\alpha\wedge\pi_X^*(\beta^{k-1})
+\binom{k}{2}\pi_S^*(\alpha^2)
 \wedge\pi_X^*(\beta^{k-2}).
$$
Using $d\beta=0$, $\partial\bar\partial\alpha=0$, and $d(\alpha^2)=0$, we obtain $\partial\bar\partial\omega^k=0.$
This also holds for $k=1$. Taking $k=n-2$ shows that the metric is astheno-Kähler for every $n\geq 3$.

We now verify directly that the product metric is not torsion-critical. For $s>0$, set
$$
\omega_s=s\,\pi_S^*\alpha+\pi_X^*\beta,
$$
and let $g_s$ be the corresponding Hermitian metric.
Explicitly, $g_s=s\,\pi_S^*g_S+\pi_X^*g_X$, where $g_X$ is the flat Kähler metric associated with $\beta$.
Thus $g_1$ is the original product metric, with fundamental
form $\omega_1=\omega$.

Since $S$ has complex dimension two,
$$
dv_{g_s}=s^2\,dv_{g_1},
\qquad
V(g_s)=s^2V(g_1).
$$
The Chern torsion comes entirely from the surface factor $S$.
Scaling $g_S$ by $s$ leaves its Chern connection unchanged
and multiplies the squared torsion norm by $s^{-1}$. Thus
$$
|T(g_s)|_{g_s}^2=s^{-1}|T(g_1)|_{g_1}^2.
$$
It follows that
$$
\int_M |T(g_s)|_{g_s}^2\,dv_{g_s}
=s\int_M |T(g_1)|_{g_1}^2\,dv_{g_1}.
$$
By the definition of $\mathcal F$, we obtain
$$
\mathcal F(g_s)
=s^{\frac{2-n}{n}}\mathcal F(g_1),
\qquad
\left.\frac{d}{ds}\right|_{s=1}\mathcal F(g_s)
=\frac{2-n}{n}\mathcal F(g_1)<0.
$$
Here $\mathcal F(g_1)>0$ because $g_1$ is not Kähler, and $n\geq 3$. Thus $g_1$ is not torsion-critical. This example shows that the torsion-critical assumption in Theorem \ref{thm:main} cannot be omitted.

\end{example}

\begin{remark}
The product calculation works more generally for any compact pluriclosed Hermitian surface and any compact K\"ahler factor. Example~\ref{ex:kodaira} shows that pluriclosedness and the astheno-K\"ahler condition together do not imply K\"ahlerness. It establishes the need for torsion-criticality in Theorem~\ref{thm:main}. It does not establish the necessity of the astheno-K\"ahler hypothesis among pluriclosed torsion-critical metrics.
\end{remark}

We now explain the role of the astheno-K\"ahler hypothesis in the comparison of first variations.

For an integer $k\geq1$, repeated use of the product rule gives
\begin{equation}\label{eq:powers}
 \partial\pb\omega^k
 =k\partial\pb\omega\wedge\omega^{k-1}
  +k(k-1)\partial\omega\wedge\pb\omega\wedge\omega^{k-2},
\end{equation}
where the second term is absent for $k=1$. 

If $g$ is pluriclosed and $n\geq4$, equation \eqref{eq:powers} becomes
\begin{equation}\label{eq:obstruction}
 \partial\pb\omega^{n-2}
 =(n-2)(n-3)\partial\omega\wedge\pb\omega\wedge\omega^{n-4}.
\end{equation}
Consequently, at a pluriclosed torsion-critical metric,
\begin{equation}\label{eq:residual}
 D\G_g(\gamma)
 =-\frac{V^{\frac{1-n}{n} }}{(n-2)!}\int_M\gamma\wedge\ddc\omega^{n-2}
\end{equation}
for every smooth real $(1,1)$-form $\gamma$. The right-hand side is the precise term left by the comparison method. Theorem~\ref{thm:G} cannot be applied unless $D\G_g$ vanishes on all such directions. 

Set $\Psi=\sqrt{-1}\partial\bar\partial\omega^{n-2}$. Taking the real $(1,1)$-form $\gamma=*\Psi$ in \ref{eq:residual} gives
$$
D\mathcal G_g(*\Psi) =-\frac{V^{\frac{1-n}{n}}}{(n-2)!} \int_M|\Psi|^2\,dv.
$$
Thus, if $\Psi \neq 0$, this variation decreases $\mathcal G$ to first order.
This is consistent with $D\mathcal F_g=0$, which alone does not imply $D\mathcal G_g=0$.

\section{Consequences under curvature assumptions}\label{sec:additional-curvature}

The rigidity theorem above works for a fixed Hermitian metric. Without extra assumptions, it cannot tell us the biholomorphism type of the complex manifold. We can get more results if we add conditions on curvature and topology. 
These results follow from Theorem \ref{thm:main}, combined with known theorems about flat manifolds and manifolds with positive holomorphic bisectional curvature. 
Readers may refer to Boothby \cite{Boothby} for the structure of compact Chern-flat Hermitian manifolds.

\begin{corollary}\label{cor:chern-flat-torus}
Let $(M^n,g)$ be a compact connected Hermitian manifold, with $n\geq3$. Suppose that $g$ is pluriclosed, astheno-K\"ahler, and torsion-critical. If its Chern curvature vanishes,
$R^{\mathrm{c}}=0$, then there is a finite unramified holomorphic covering
\[
\mathbb C^n/\Lambda\longrightarrow M,
\]
where $\Lambda\subset\mathbb C^n$ is a lattice of real rank $2n$.
If, in addition, $b_1(M)=2n$, then $M$ itself is biholomorphic to a complex torus. Here $b_1(M)=\dim_{\mathbb R}H^1_{\mathrm{dR}}(M;\mathbb R)$.
\end{corollary}

\begin{proof}
By Theorem \ref{thm:main}, g is Kähler. Its Chern connection coincides with the Levi-Civita connection. Thus, $R^{\mathrm{c}}=0$ means g is flat.
The universal cover inherits the lifted metric. It is complete, simply connected and flat, so it is isometric to $\mathbb R^{2n}$.
The lifted complex structure is parallel. We may choose this isometry to be holomorphic, and it identifies the universal cover with $\mathbb C^n$.

Let $\Gamma$ be the deck transformation group. Each $\delta\in\Gamma$ has the form
$$
\delta(z)=A_\delta z+a_\delta,
\qquad A_\delta\in U(n),\quad a_\delta\in\mathbb C^n.
$$
By the first Bieberbach theorem, the translation subgroup of $\Gamma$ is a full lattice $\Lambda$ of finite index; see \cite{Charlap86}. Thus
$$
\mathbb C^n/\Lambda\longrightarrow\mathbb C^n/\Gamma=M
$$
is a finite covering. The map is holomorphic and locally biholomorphic, so it is unramiﬁed.

Suppose now that $b_1(M)=2n$. By Hodge theory, the space of harmonic real $1$-forms has dimension $2n$. 

For a real 1-form $\alpha$, by the Weitzenb\"ock--Bochner formula in \cite{Petersen},
$$
\Delta_d \alpha = (\nabla^{\text{LC}})^*\nabla^{\text{LC}} \alpha + \mathrm{Ric}\cdot \alpha,
$$
where $\Delta_d = dd^* + d^*d.$ If $\alpha$ is harmonic, then $\Delta_d \alpha = 0$. Take the inner product with $\alpha$ and integrate:
$$
0 = \int_M \langle (\nabla^{\text{LC}})^*\nabla^{\text{LC}} \alpha,\alpha\rangle dv + \int_M \left\langle \mathrm{Ric}\cdot \alpha ,\alpha \right\rangle \,dv.
$$
Integration by parts gives
$$
\int_M \langle (\nabla^{\text{LC}})^*\nabla^{\text{LC}} \alpha,\alpha\rangle\,dv = \int_M |\nabla^{\text{LC}} \alpha|^2\,dv.
$$
Since the metric is flat, the Ricci curvature vanishes, so
$$
0 = \int_M |\nabla^{\text{LC}} \alpha|^2\,dv.
$$
The integrand is non-negative and continuous, hence $\nabla^{\text{LC}} \alpha = 0.$

Thus every harmonic $1$-form is parallel. Evaluation at any point is injective on the space of parallel one-forms: a parallel form that vanishes at one point vanishes everywhere on the connected manifold.  
Hence the $2n$ linearly independent harmonic $1$-forms form a real coframe at every point.

Their lifts to $\mathbb C^n$ are constant $1$-forms invariant under all deck transformations. 
Thus $A_\delta^*$ fixes a basis of the real cotangent space, so $A_\delta=I$ for every $\delta\in\Gamma$. Every deck transformation is therefore a translation. 
Hence $\Gamma=\Lambda$,and $M=\mathbb C^n/\Lambda$ is a complex torus.
\end{proof}

\begin{corollary}\label{cor:positive-bisectional-projective}
Let $(M^n,g)$ be a compact connected Hermitian manifold, with $n\geq3$. Suppose that $g$ is pluriclosed, astheno-K\"ahler, and torsion-critical. If its Chern holomorphic bisectional curvature is strictly positive, then $M$ is biholomorphic to $\mathbb{CP}^n$.
\end{corollary}

\begin{proof}
Theorem~\ref{thm:main} implies that $g$ is K\"ahler. In particular, its Chern connection agrees with the connection induced by the Levi-Civita connection, so the Chern holomorphic bisectional curvature is precisely the usual K\"ahler holomorphic bisectional curvature.

For any point $p\in M$, define the bisectional curvature
$$
B_g(X,Y)=\frac{R_g(X,\overline{X},Y,\overline{Y})}{|X|_g^2|Y|_g^2},\qquad X,Y\neq 0,\qquad X,Y\in T_p^{1,0}M.
$$
Since the denominator is strictly positive,
$$
B_g(X,Y)>0 \iff R_g(X,\overline{X},Y,\overline{Y})>0.
$$

The Siu--Yau theorem \cite{SiuYau80} states that a compact K\"ahler
manifold with strictly positive holomorphic bisectional curvature is
biholomorphic to complex projective space. Applying that theorem gives
$M\cong\mathbb{CP}^n$.
\end{proof}

\begin{remark}\label{rem:classification-scope}
  Corollary \ref{cor:chern-flat-torus} gives a finite covering of $M$ by a complex torus. It does not imply that $M$ itself is a complex torus.
  Indeed, the deck transformations of the universal cover may have nontrivial linear parts. The condition $b_1(M)=2n$ forces all these linear parts to be the identity.
  Then $M$ is a complex torus.
  
  In this Chern-flat case, the astheno-K\"ahler and torsion-critical assumptions can be omitted. Compact Chern-flat Hermitian metrics are balanced; see \cite{ZZ25} and the references therein.
  Thus pluriclosedness and Corollary \ref{cor:balanced} already imply that $g$ is K\"ahler.
\end{remark}

\noindent\textbf{Acknowledgments.} We would like to express our sincere gratitude to Professor Fangyang Zheng for his insightful guidance, valuable discussions, and continued support throughout the course of this research.

\noindent\textbf{Disclosure on AI assistance.} The research problem and the principal mathematical ideas were formulated by the authors. AI tools were used only during the final stages of manuscript preparation to improve the exposition.The authors retain full responsibility for the correctness of the final work.

\end{document}